\documentclass{article}

\usepackage{a4wide} 
\usepackage{amsmath} 
\usepackage{amsfonts} 
\usepackage{epsfig} 
\usepackage{latexsym} 
\usepackage{psfrag} 
\usepackage{graphicx}

\newtheorem{lemma}{Lemma} 
\newtheorem{proposition}{Proposition}

\newcommand{\eps}{\varepsilon}   
\newcommand{\qed}{\rightline{$\Box$}} 

\newcommand{\one}{{\bf 1}}

\newcommand{\onesccn}{{\bf 1}}

\def\R{{\mathbb R}}

\def\P{{\mathbb P}}

\def\u{{u_T}}
\def\v{{v}}

\begin{document} 
 
\title{A singular perturbation approach\\ for choosing PageRank damping factor}

\author{
Konstantin Avrachenkov\thanks{INRIA Sophia Antipolis, France,
E-mail: k.avrachenkov@sophia.inria.fr}, 
Nelly Litvak\thanks{University of Twente, The Netherlands,
E-mail: n.litvak@utwente.nl} and 
Kim Son Pham\thanks{
St.Petersburg State University, Russia,
E-mail: sonsecure@yahoo.com.sg}
}

\date{}

\maketitle
 
\begin{abstract}
The choice of the PageRank damping factor is not evident. The Google's choice for the value $c=0.85$ 
was a compromise between the true reflection of the Web structure and numerical efficiency. 
However, the Markov random walk on the original Web Graph does not reflect the importance 
of the pages because it absorbs in dead ends. Thus, the damping factor is needed not only for 
speeding up the computations but also for establishing a fair ranking of pages. 
In this paper, we propose new criteria for choosing the damping factor, based on the 
ergodic structure of the Web Graph and probability flows. Specifically, we require that 
the core component receives a fair share of the PageRank mass. Using singular perturbation approach
we conclude that the value $c=0.85$ is too high and suggest that the damping factor should be
chosen around 1/2. 
As a by-product, we describe the ergodic structure of the OUT component of the Web Graph in
detail. 
Our analytical results are confirmed by experiments on two large samples of the Web Graph.

\bigskip

\end{abstract}


\section{Introduction}

Surfers on the Internet frequently use search engines to find pages satisfying their query.
However, there are typically hundreds or thousands of relevant pages available on the
Web. Thus, listing them in a proper order is a crucial and non-trivial task. One can
use several criteria to sort relevant answers. It turns out that the link-based criteria
 provide rankings that appear to be
very satisfactory to Internet users. The examples of link-based criteria are PageRank~\cite{PB98}
used by search engine Google, HITS \cite{K99} used by search engines Teoma and Ask, and
SALSA \cite{LM00}. 

In the present work we restrict ourselves to the analysis of the PageRank
criterion and use the following definition of PageRank from~\cite{LM06}. Denote by $n$ the total
number of pages on the Web and define the $n\times n$ hyperlink matrix $P$ as follows:
\begin{equation}
\label{P}
p_{ij} =
\left\{ \begin{array}{ll}
1/d_i, & \mbox{if page $i$ links to $j$},\\
1/n, & \mbox{if page $i$ is dangling},\\ 
0, & \mbox{otherwise},
\end{array} \right.
\end{equation}
for $i,j=1,...,n$, where $d_i$ is the number of outgoing links from page $i$. We recall that the page is
called dangling if it does not have outgoing links.
In order to make the hyperlink graph connected, it is assumed that at each step, with some probability, a  random surfer goes to an arbitrary Web page sampled from the uniform distribution. Thus, the PageRank is defined as a stationary
distribution of a Markov chain whose state space is the set of all Web pages, and the
transition matrix is
\begin{equation}
\label{GoogleMatrix}
G = cP + (1-c)(1/n)E, 
\end{equation}
where $E$ is a matrix whose all entries are equal to one, and $c \in (0,1)$ is a probability of following a hyperlink. The constant $c$ is often 
referred to as a damping factor. The Google matrix $G$ is stochastic, aperiodic, and irreducible, 
so there exists a unique row vector $\pi$ such that
\begin{equation}
\label{BalanceEq}
\pi G = \pi, \quad \pi \one =1,
\end{equation}
where $\one$ is a column vector of ones. The row vector $\pi$ satisfying (\ref{BalanceEq}) 
is called a PageRank vector, or simply PageRank. If we consider a surfer that follows a hyperlink with probability $c$ and jumps to a random page with probability $1-c$, then $\pi_i$ can be interpreted as 
a stationary probability that the surfer is at page $i$.

The damping factor $c$ is a crucial parameter in the PageRank definition. It regulates
the level of the uniform noise introduced to the system. Based on the publicly available
information Google originally used $c=0.85$. There is the following empirical explanation of this
choice, see e.g. \cite{LM06}: it seems that the closer the value of the damping factor to 
one, the better the graph structure of the Web is represented in the PageRank vector.
However, when the value of $c$ approaches one, the rate of power iteration method slows
down significantly. The choice $c=0.85$ appears to be a reasonable compromise between
the two antagonistic objectives.  However, in \cite{Boldi1} the authors argue that choosing
the value of $c$ too close to one is not necessarily a good thing to do. Not only the
power iteration method becomes very slowly convergent but also the ranking of the important
pages becomes distorted. Independently of \cite{Boldi1}, this phenomenon was also mentioned in~\cite{AL06}. We also remark that another incentive to reduce $c$ is that it will increase the robustness of the PageRank towards small changes in the link structure. That is, with smaller $c$, one can bound the influence of outgoing links of a page (or a small group of pages) on the PageRank of other groups~\cite{Bianchini} and on its own PageRank~\cite{AL06}.

In the present work, we go further than \cite{Boldi1, AL06} and suggest that even the value
$c=0.85$ is by far too large. Our argument is that one has to make a choice of $c$ to reflect
the natural intensity of the probability flow in the absorbing Markov chain associated with the
Web Graph. Our argument is based on the singular perturbation theory \cite{A99,KT93,PG88,YZ05}. 
It turns out that
the value of $c$ that adequately reflects the flow of probability is very close to $1/2$. We note that the value $c=1/2$ was used in \cite{PRcitations} to find gems in scientific citations, where the authors justified this choice by intuitive argument discussed in more detail in Section~\ref{sec:conclusions}. In this work, we present a mathematical evidence for setting $c=1/2$ in the PageRank formula.

Of course, a drastic reduction of $c$ considerably accelerates the computation
of PageRank by numerical methods \cite{ALNO06,B05,LM06}. We would like to mention
that choosing smaller value for the damping factor could have similar effect on numerical
methods as choosing fast decreasing damping function \cite{Boldi2}. 

As a by-product of the application of the singular perturbation approach we obtain 
a refinement of the graph structure of the Web. We demonstrate that the dead-end strongly 
connected components have unjustifiably large PageRank with damping factor $c=0.85$ and
by taking $c=0.5$ one can mitigate this problem.  The results presented in this work are confirmed by experimental data that we obtained from two large samples of the Web Graph, described in Section~\ref{sec:datasets}.

The main contributions of this paper are as follows. First, in Section~\ref{sec:ergodic}, we describe the ergodic structure of the Web Graph and show how this structure changes under assumption that the dangling pages have a link to all pages in the Web, as in (\ref{P}). In particular, we discover an Extended Strongly Connected Component (ESCC) that contains a majority of the Web pages. Using the theory of singular perturbations, we find an exact formula for the limiting PageRank distribution when $c\to 1$. This result immediately implies that the limiting PageRank mass of ESCC equals zero. Next, in Section~\ref{sec:bounds}, we analytically characterize the PageRank mass of ESCC as a function of $c$, and we obtain simple bounds for this function. Further, in Section~\ref{sec:c=1/2}, we argue that $c=1/2$ ensures that ESCC receives a fair share of total PageRank mass. We conclude with a short discussion of the present results and future research directions in Section~\ref{sec:conclusions}.

\section{Datasets}
\label{sec:datasets}

For our numerical experiments, we have collected two Web Graphs, which we denote by INRIA and FrMathInfo. 
The Web Graph INRIA was taken from the site of INRIA, the French Research Institute of Informatics and
Automatics. The seed for the INRIA collection was Web page {\tt www.inria.fr}. It is a typical large Web
site with around 300.000 pages and 2 millions hyperlinks. We have crawled the INRIA site until we have 
collected all pages belonging to INRIA.

The Web Graph FrMathInfo was crawled with 
the initial seeds of 50 mathematics and informatics laboratories of France, taken from Google Directory. The crawl was executed by breadth first search and 
the depth of this crawl was 6. The FrMathInfo Web Graph contains around 700.000 pages and  
8 millions hyperlinks. We expect our datasets to be enough representative. This is justified by 
the fractal structure of the Web \cite{self-similar}.

The link structure of these two Web Graphs is stored in Oracle database. Due to sparsity
of the Web Graph and reasonable sizes of our datasets, we can store the adjacency lists in RAM to speed up 
the computation of PageRank and other quantities of interest. This enables us to make 
more iterations, which is extremely important in the case when the damping factor $c$ is close to one. 
Our PageRank computation program consumes about one hour to make 500
iterations for the FrMathInfo dataset and about haft an hour for the INRIA dataset for the same
number of iterations. Our algorithms for discovering the ergodic structures of the Web Graph are based on 
Breadth First Search and Depth First Search methods, which are
linear in the sum of number of nodes and links.

\section{Ergodic structure of the Web Graph}
\label{sec:ergodic}

In \cite{WebGraph1,WebGraph2} the authors have studied the graph structure of the Web. 
In particular, in \cite{WebGraph1,WebGraph2} 
it was shown that the Web Graph can be divided into three principle components: the Giant Strongly
Connected Component, to which we simply refer as SCC component, the IN component and the OUT component. 
The SCC component is the largest strongly connected component in the Web Graph. In fact, it is 
larger than the second largest strongly connected component by several orders of magnitude.
Following hyperlinks one can come from the IN component to the SCC component but it is not possible to return
back. Then, from the SCC component one can come to the OUT component and it is not possible to
return to SCC from the OUT component. 

With this structure in mind, we would like to analyze ergodic properties of the random walk on the Web Graph. If a node has outgoing links, then such random walk follows one of these links with uniform distribution. However, as in the definition of PageRank, we have to define how the process evolves when it reaches one of dangling nodes. This choice has a crucial influence on the ergodic structure of the associated Markov chain. There are three natural possibilities: 1)~the process is absorbed in the
dangling node; 2)~the process moves to the predecessor node, or 3)~the process moves to an
arbitrary node. In this paper we focus on the latter option, which is used in the original PageRank model \cite{PB98}. The first two options are definitely worthy to be considered as well, and it is a nice topic for future research. 

\begin{figure}[hbt]
              \centering {\epsfxsize=2.8in \epsfbox{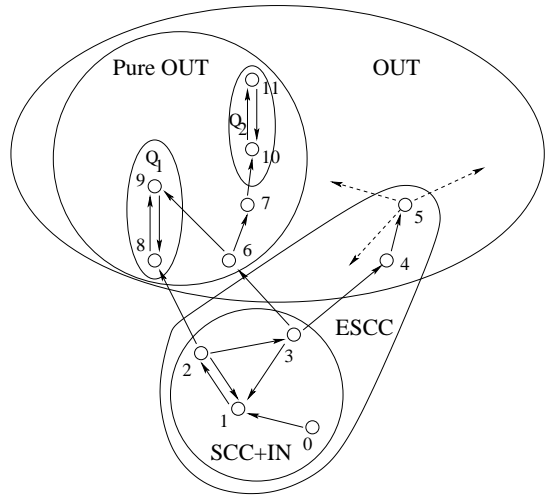}}
                \caption{\small Example of a graph}
        \label{fig:smallgraph}
\end{figure}

Thus, throughout the paper we consider a random walk with transition matrix $P$ given by (\ref{P}). As we shall see below, the analysis of the ergodic structure of $P$ leads to a more detailed description of the OUT component, and it allows us to evaluate the effect of damping factor
on PageRank. Obviously, the graph induced by $P$ has a much higher connectivity than the original Web Graph. In particular, if the random walk
can move from a dangling node to an arbitrary node with the uniform distribution, then the 
Giant SCC component increases further in size. We refer to this new strongly
connected component as the Extended Strongly Connected Component (ESCC). First, we note that due to the artificial links from the dangling nodes,
the SCC component and IN component are now inter-connected and are parts of the Extended
SCC. Then, if there are dangling nodes emanating from some nodes in the OUT component,
these nodes together with all their predecessors become a part of the Extended SCC. 
Let us consider an example of the graph presented in Figure~\ref{fig:smallgraph}.
Node~0 represents the IN component, nodes from 1 to 3 form the SCC component, and
the rest of the nodes, nodes from 4 to 11, are in the OUT component. Node~5 is a dangling node, thus,
artificial links go from the dangling node~5 to all other nodes. After addition of the
artificial links, all nodes from 0 to 5 form a new strongly connected component, which is the ESCC
in this example.

By renumbering the nodes, the transition matrix $P$ can be then transformed to the
following form   
\begin{equation}
\label{ESCC}
P=\left[ \begin{array}{cc} 
Q & 0  \\
R & T \end{array} \right],
\end{equation}
where the block $T$ corresponds to the Extended SCC, the block $Q$ corresponds to
the part of the OUT component without dangling nodes and their predecessors,  and the block $R$ corresponds to the transitions from ESCC to the nodes in block $Q$. We refer to the set of nodes in the block $Q$ as Pure OUT component.
In the example of graph on Figure~\ref{fig:smallgraph} the Pure OUT component consists of
nodes from 6 to 11. Typically, the Pure OUT component is much smaller than the Extended SCC
component. The sizes of all components for our two datasets are displayed in Table~1.
We would like to note that the zero size of the IN components should not come as a surprise.
To crawl the Web Graph we have used the Breadth First Search method and have started from
important pages. Therefore, it is natural that the seed pages belong to the Giant SCC and
there is no IN component. For the purposes of the present research the absence of the IN
component is not a problem as the dangling nodes unite the IN and the Giant SCC into the
Extended SCC.  

As was observed in \cite{moler}, the PageRank vector can be expressed by the following
formula
\begin{equation}
\label{PRformula}
\pi=\frac{1-c}{n} \one^T [I-cP]^{-1}.
\end{equation}
If we substitute the expression (\ref{ESCC}) for the transition matrix $P$
into (\ref{PRformula}), we obtain the following formula for the part of the 
PageRank vector corresponding to the nodes in ESCC:
$$
\pi_{T}=\frac{1-c}{n} \one^T [I-cT]^{-1}, 
$$ 
or, equivalently,
\begin{equation}
\label{PRESCC}
\pi_{T} = \alpha (1-c) u_{T} [I-cT]^{-1}, 
\end{equation} 
where $\alpha=n_T/n$ and $n_T$ is the number of nodes in ESCC, 
and where $u_T$ is the uniform distribution over all ESCC nodes.
We shall also use $n_Q=n-n_T$, which is the number of nodes in the 
Pure OUT component.

First, we note that since matrix $T$ is substochastic, the inverse
$[I-T]^{-1}$ exists and consequently $\pi_T \to 0$ as $c \to 1$.
Clearly, as was also observed in \cite{Boldi1}, it is not good to take
the value of $c$ too close to one. Next, we argue that even the value
of $0.85$ is too large.

\begin{figure}[hbt]
              \centering {\epsfxsize=3.0in \epsfbox{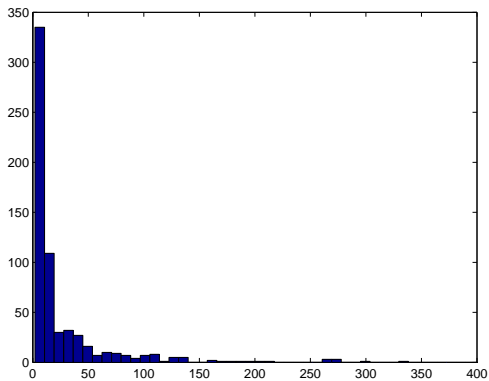}}
                \caption{\small Histogram of SCCs' sizes of Pure OUT, INRIA dataset}
        \label{fig:histINRIA}
\end{figure}

\begin{figure}[hbt]
              \centering {\epsfxsize=3.0in \epsfbox{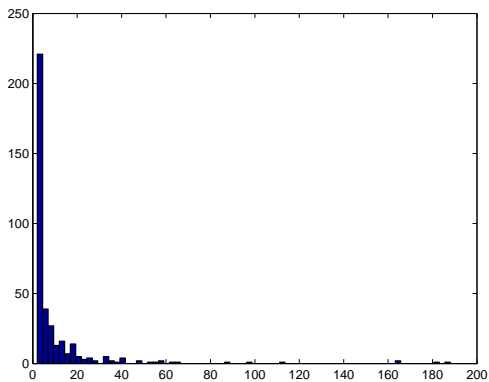}}
                \caption{\small Histogram of SCCs' sizes of Pure OUT, FrMathInfo dataset}
        \label{fig:histFrMathInfo}
\end{figure}

Let us analyze the structure of the Pure OUT component in more detail.
It turns out that there are many disjoint strongly connected components
inside the Pure OUT component. One can see the histograms of the SCCs'
sizes of the Pure OUT for two our datasets INRIA and FrMathInfo in 
Figures~\ref{fig:histINRIA}~and~\ref{fig:histFrMathInfo}. In particular,
there are many SCCs of size 2 and 3 in the Pure OUT component.

By appropriate renumbering of the states, we can refine (\ref{ESCC}) as
follows:
\begin{equation}
\label{ESCCdetail0}
P=\left[ \begin{array}{ccccc} 
Q_1 &        &     &     & 0 \\
    & \ddots &     &     &  \\
0   &        & Q_m &     &  \\
S_1 & \cdots & S_m & S_0 & 0 \\    
R_1 & \cdots & R_m & R_0 & T 
\end{array} \right],
\end{equation}

For instance, in example of the graph from Figure~\ref{fig:smallgraph}, the nodes
8 and 9 correspond to block $Q_1$, nodes 10 and 11 correspond to block $Q_2$, and
nodes 6 and 7 correspond to blocks $S$. 

Since the random walk will be eventually absorbed in one of the $Q$ blocks,
we can simplify notations for our further analysis. Namely, define the submatrices
$$
\tilde{R}_i=\left[\begin{array}{c}S_i\\R_i\end{array}\right],\;i=1,\ldots,m; \quad \tilde{T}=
\left[\begin{array}{cc}S_0&0\\R_0&T\end{array}\right].
$$
Then the structure (\ref{ESCCdetail0}) becomes 
 \begin{equation}
\label{ESCCdetail}
P=\left[ \begin{array}{cccc} 
Q_1 &        & 0   & 0 \\
    & \ddots &     &   \\
0   &        & Q_m & 0 \\
\tilde{R}_1 & \cdots & \tilde{R}_m & \tilde{T} 
\end{array} \right].
\end{equation}


Next, we note
that if $c<1$, then the Markov chain induced by matrix $G$ is ergodic.
However, if $c=1$, the Markov chain becomes non-ergodic. In particular,
if the process moves to one of the $Q_i$ blocks,
it will
never leave this 
block. Hence, the random walk governed by the Google transition
matrix (\ref{GoogleMatrix}) is in fact a singularly perturbed Markov
chain. 

According to the singular perturbation theory (see e.g., \cite{A99,KT93,PG88,YZ05}),
the PageRank vector goes to
some limit as the damping factor goes to one. Using the results 
of the singular perturbation theory we can characterize explicitly
this limit.

\begin{proposition}
Let $\mu_i$ be a limiting stationary distribution of the Markov process
governed by 
$P$ when the process settles in $Q_i$, the $i$-th SCC of 
the Pure OUT component. Namely, vector $\mu_i$ is a unique solution
of the equations
$$
\mu_i Q_i = \mu_i, \quad \mu_i \one = 1.
$$
Then, we have
$$
\lim_{c \to 1} \pi(c)=
\left[\bar{\pi}_1 \ \cdots \ \bar{\pi}_m \ 0 \right],
$$
where
\begin{equation}
\label{barpi}
\bar{\pi}_i=
\left(\frac{n_i}{n}+\frac{n_{\tilde{T}}}{n}u_{\tilde{T}}[I-\tilde{T}]^{-1}\tilde{R}_i\one\right)\mu_i,
\end{equation}
for $i=1,...,m$ and the zeros at the end of the limiting vector correspond to all nodes,
which are not in $Q_i$, $i=1,\ldots,n$, that is, not in any SCC of the Pure OUT component. 
\end{proposition}
{\bf Proof:}
First, we note that if we make a change of variables $\eps = 1-c$ the Google matrix
becomes a transition matrix of a singularly perturbed Markov chain as in Lemma~\ref{lm:SPMC}
with $C=\frac{1}{n}\one \one^T -P$. Let us calculate the aggregated generator matrix $D$.
$$
D=MCQ=\frac{1}{n}\one\one^TQ-MPQ
$$
Using $MP=M$, $MQ=I$, and $M\one=\one$ where vectors $\one$ are of appropriate dimensions,
we obtain
$$
D=\frac{1}{n}\one\one^TQ-I=
$$
$$
\frac{1}{n}\one[n_1+n_{\tilde{T}}u_{\tilde{T}}[I-\tilde{T}]^{-1}\tilde{R}_1\one,\cdots,
n_m+n_{\tilde{T}}u_{\tilde{T}}[I-\tilde{T}]^{-1}\tilde{R}_m\one]-I
$$
Since the aggregated transition matrix $D+I$ has identical rows, its stationary
distribution $\nu$ is just equal to these rows. Thus, invoking Lemma~\ref{lm:SPMC} 
we obtain (\ref{barpi}).   
\qed

The second term inside the brackets in formula (\ref{barpi}) corresponds to the PageRank
mass that an SCC component in Pure OUT receives from the Extended SCC. If $c$ is close to
one, then this contribution can outweight by far the fair share of the PageRank
which is given by $\frac{n_i}{n}$. For instance, in our numerical experiments with $c=0.85$, the PageRank mass of the Pure OUT component in the INRIA dataset equals $1.95n_Q/n$, whereas a `fair share' is $n_Q/n$. In the other dataset, FrMathInfo, the unfairness  is even more pronounced: the PageRank mass of the Pure OUT component is $3.44n_Q/n$.  This gives users an incentive to create `dead-ends': groups of pages that link only to each other. In the next sections we quantify the influence of parameter $c$ and show that in order to obtain balanced probability flow, $c$ should be taken around 1/2.

\begin{table}[htb]
\label{tab:table1}
\centerline{\begin{tabular}{|r|r|r|}
\hline
$ $&$INRIA$&$FrMathInfo$\\
\hline \hline
Total size & 318585 & 764119\\
Number of nodes in SCC  & 154142 & 333175\\
Number of nodes in IN  & 0 & 0 \\
Number of nodes in OUT  & 164443 & 430944 \\
Number of nodes in ESCC & 300682 & 760016\\
Number of nodes in Pure OUT & 17903 & 4103\\
Number of SCCs in OUT & 1148 & 1382\\
Number of SCCs in Pure Out & 631 & 379\\
\hline
\end{tabular}}
\caption{Component sizes in INRIA and FrMathInfo datasets}
\end{table}

\section{PageRank mass of ESCC}
\label{sec:bounds}
 
Let us consider the PageRank mass of the Extended SCC component (ESCC) described in the previous section. Thus, we continue to analyze the transition matrix in the form  presented in (\ref{ESCC}). 
Our goal now is to characterize the behavior of the total PageRank mass of the ESCC component as a function of $c \in [0,1]$. From (\ref{PRESCC}) we have
\begin{align}
\nonumber ||\pi_T(c)||_1&=
\pi_{T}(c)\onesccn=(1-c)\alpha\u[I-cT]^{-1}\onesccn\\
\label{piscc_sum}
&=(1-c)\alpha\u\sum_{k=0}^\infty c^kT^k\onesccn.\end{align} 
Clearly, since $T$ is substochastic, we have $||\pi_{T}(0)||_1=\alpha$ and $||\pi_{T}(1)||_1=0$. Also, it is easy  to show that 
\[\frac{d}{dc}||\pi_{T}(c)||_1=-\alpha\u[I-cT]^{-2}[I-T]\onesccn<0\]
and
\[\frac{d^2}{dc^2}||\pi_{T}(c)||_1=-2\alpha\u[I-cT]^{-3}T[I-T]\onesccn<0.\]
Hence, $||\pi_{T}(c)||_1$ is a concave decreasing function. 

In order to get a better idea about the behavior of this function, we derive a series of bounds. If we define\[\underline{p}=\inf_{k\ge 1}[uT^k\onesccn]^{1/k},\quad 
\overline{p}=\sup_{k\ge 1}[uT^k\onesccn]^{1/k},\]
then it follows immediately from (\ref{piscc_sum}) that
\begin{equation}
\label{b1}
\frac{\alpha(1-c)}{1-c\underline{p}}\le ||\pi_T(c)||_1\le\frac{\alpha(1-c)}{1-c\overline{p}}.
\end{equation}
Now, let $\lambda_1$ be the Perron-Frobenius eigenvalue of $T$, and let $\tau$ be a random time when a random walk induced by $T$ leaves  ESCC given that the  initial distribution is uniform on ESCC. It is well known that 
\[\lambda_1=\lim_{k\to\infty}\P[\tau>k|\tau>k-1]=\lim_{k\to\infty}\frac{uT^k\onesccn}{uT^{k-1}\onesccn}\]
Thus, we evaluate $\lambda_1$ iteratively by computing
\begin{equation}
\label{lambdak}
\lambda_1^{(k)}=\frac{uT^k\onesccn}{uT^{k-1}\onesccn},\quad k\ge 1,\end{equation}
where the numerator and denominator are simply results of the power iterations of $T$.
From the definition of $\lambda_1^{(k)}$ it is easy to see that if the sequence $\lambda_1^{(k)}$, $k\ge 1$, is increasing then
the sequence $(uT^k\onesccn)^{1/k}$, $k\ge 1$, is also increasing, and thus in this case $\overline{p}=\lambda_1$ and $\underline{p}=p_1$, where $p_1=\u T\onesccn=\P(\tau>1)$. Then equation (\ref{b1}) becomes
\begin{equation}
\label{b}
\frac{\alpha(1-c)}{1-cp_1}\le ||\pi_T(c)||_1\le\frac{\alpha(1-c)}{1-c\lambda_1}.\end{equation}  
Although in our experiments we indeed observed that the sequence $\lambda_1^{(k)}$, $k\ge 1$, is increasing for both INRIA and FrMathInfo datasets, this condition is still too strong and we presume that it may fail in some cases. In the next proposition we provide much milder and more intuitive conditions under which (\ref{b}) still holds. 

\begin{proposition} 
\label{prop1}
Let $\lambda_1$ be the Perron-Frobenius eigenvalue of $T$, and define
$p_1=\u T\onesccn$. 
\begin{itemize}
\item[(i)] 
If $p_1<\lambda_1$ then  
\begin{equation}
\label{ub}
||\pi_{T}(c)||_1<\frac{\alpha(1-c)}{1-c\lambda_1},\quad c\in(0,1).\end{equation}
\item[(ii)] If $1/(1-p_1)<\u[I-T]^{-1}\onesccn$ then  
\begin{equation}
\label{lb}
||\pi_{T}(c)||_1>\frac{\alpha(1-c)}{1-cp_1},\quad c\in(0,1).\end{equation}
\end{itemize}
\end{proposition}

{\bf Proof.} (i) The function $f(c)=\alpha(1-c)/(1-\lambda_1 c)$ is decreasing and concave, and so is $||\pi_{T}(c)||_1$. Also, $||\pi_{T}(0)||_1=f(0)=\alpha$, and $||\pi_{T}(1)||_1=f(1)=0$. Thus, for $c\in (0,1)$, the plot of $||\pi_{T}(c)||_1$ is either entirely above or entirely below $f(c)$. In particular, if the first derivatives satisfy $||\pi'_{T}(0)||_1<f'(0)$, then $||\pi_{T}(c)||_1<f(c)$ for any $c\in(0,1)$. Since $f'(0)=\alpha(\lambda_1-1)$ and $||\pi'_{T}(0)||_1=\alpha(p_1-1)$, we see that 
 $p_1<\lambda_1$ implies (\ref{ub}). 

The proof of (ii) is similar. We consider a concave decreasing function $g(c)=\alpha(1-c)/(1-p_1 c)$ and note that $g(0)=\alpha$, $g(1)=0$. Now, if the condition in (ii) holds then $g'(1)>||\pi'_{T}(1)||_1$, which implies (\ref{lb}). \qed

Note that the conditions of Proposition~\ref{prop1} are satisfied when the sequence $\lambda_1^{(k)}$, $k\ge 1$, is increasing in $k$. 

The condition $p_1<\lambda_1$, which gives the upper bound, has a clear intuitive interpretation. Let $\tilde{\pi}_T$ be a quasi-stationary distribution of $T$. By definition, $\tilde{\pi}_T$ is the probability-normed left Perron-Frobenius eigenvector of $T$, and it is well-known that $\tilde{\pi}_T$ is a limiting probability distribution obtained under condition that the random walk does not leave the ESCC component (see e.g. \cite{seneta}). Hence, $\tilde{\pi}_TT=\lambda_1\tilde{\pi}_T$, and the condition $p_1<\lambda_1$ means that the chance to stay in ESCC for one step in the quasi-stationary regime is higher than starting from the uniform distribution $\u$.  This inequality looks quite natural, since the quasi-stationary distribution should somehow favor states, from which the chance to leave ESCC is lower. Therefore, although $p_1<\lambda_1$ does not hold in general, one may expect that it should hold for transition matrices describing large entangled graphs.

\begin{figure}[hbt]
              \centering {\epsfxsize=2.8in \epsfbox{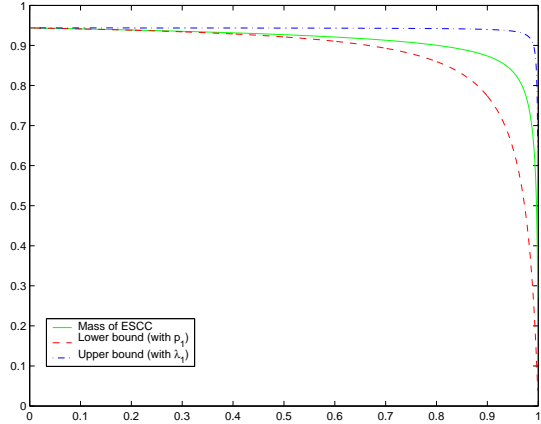}}
                \caption{\small PageRank mass of ESCC and bounds, INRIA dataset}
        \label{fig:INRIA_escc}
\end{figure}

\begin{figure}[hbt]
              \centering {\epsfxsize=2.8in \epsfbox{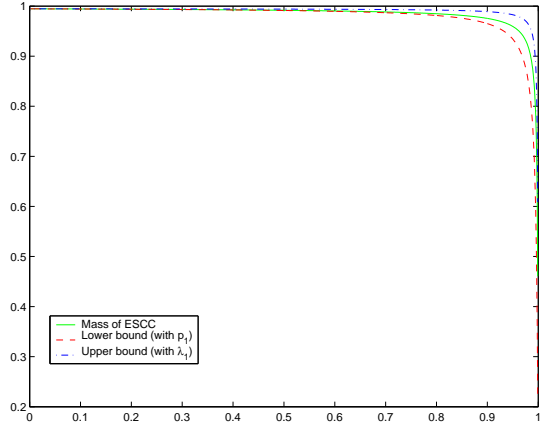}}
                \caption{\small PageRank mass of ESCC and bounds, FrMathInfo dataset}
        \label{fig:FrMathInfo_escc}
\end{figure}

With the help of the derived bounds we can conclude that the function $||\pi_{T}(c)||_1$  decreases very slowly for small and moderate values of $c$, and it decreases extremely fast when $c$ becomes close to 1. This typical behavior is clearly seen in Figures~\ref{fig:INRIA_escc},~\ref{fig:FrMathInfo_escc}, where $||\pi_{T}(c)||_1$ is plotted with a solid line. 
In order to evaluate $||\pi_{T}(c)||_1$ we did not compute it separately for different values of $c$ but rather presented it as a function of $c$ so that {\it any} value of $c$ could be substituted. For that, we stored the values $\lambda_1^{(k)}$, $k\ge 1$, and then used (\ref{piscc_sum}) and (\ref{lambdak}) to obtain
\[||\pi_{T}(c)||_1=\alpha\sum_{k=0}^{\infty}c^k\prod_{l=1}^k\lambda_1^{(k)}, \quad c\in[0,1].\]
This is a more direct approach compared to \cite{Boldi1}, where the authors used derivatives of the PageRank to present $\pi$ as a function of $c$.
As for the bounds, the values $\lambda_1$ and $p_1$ can be directly substituted in (\ref{ub}) and (\ref{lb}), respectively. For the INRIA dataset we have $p_1=\lambda_1^{(1)}=0.97557$, $\lambda_1=0.99954$, and for the FrMathInfo dataset we have $p_1=0.99659$, $\lambda_1=0.99937$.  

In the next section we use the above results on $||\pi_{T}(c)||_1$ and its bounds to determine the values of $c$ that reflect natural probability flows through the ESCC component.

\section{Why the damping factor should be 1/2}
\label{sec:c=1/2}

Since ESCC is by far more important and interesting part of the Web than the Pure~OUT component, it would be reasonable to ensure that the PageRank mass of ESCC is at least the fraction of nodes in this component (we denoted this fraction by $\alpha$). However, because $||\pi_{T}(c)||_1$ is decreasing, and $||\pi_{T}(0)||_1=\alpha$, 
it follows that the total PageRank mass of  ESCC is smaller than $\alpha$ for any value $c>0$. 

Now let us discuss an `optimal' choice of $c$. First of all, $c$ can not be too close to one because in this case the PageRank mass of the giant ESCC component will be close to 0. This was observed independently in \cite{AL06,Boldi1}. Specifically, from the analysis above it follows that the value of $c$ should not be chosen in the critical region where the PageRank mass of the ESCC component is rapidly decreasing. Luckily, the shape of the function $||\pi_{T}(c)||_1$ is such that it decreases drastically only when $c$ is really close to one, which leaves a lot of freedom for choosing $c$. In particular, the famous Google constant $c=0.85$ is small enough to ensure a reasonably large PageRank mass of ESCC. 

However, as we have observed in Section~\ref{sec:ergodic}, even moderately large values of $c$ result in an unfairly large PageRank mass of the Pure OUT component. Now, our goal is to find the values of $c$ that lead to a `fair' distribution of the PageRank mass between the Pure OUT and the ESCC components. 

Formally, we would like to define a number $\gamma\in(0,1)$ such that a desirable PageRank mass of ESCC could be written as $\gamma\alpha$, and then find the value $c^*$ that satisfies
\[||\pi_{T}(c^*)||_1=\gamma\alpha.\]
Then $c\le c^*$  will ensure that $||\pi_{T}(c)||_1\ge\gamma\alpha$. Naturally, $\gamma$ should somehow reflect the properties of the substochastic block $T$. For instance, as $T$ becomes closer to stochastic matrix, $\gamma$ should also increase. One possibility to do it is to define
\[\gamma={\v}T\onesccn,\] 
where $\v$ is a row vector representing some probability distribution on ESCC. Then the damping factor $c$ should satisfy
\[c\le c^*,\]
where $c^*$ is given by
\begin{equation}
\label{c*}
||\pi_{T}(c^*)||_1=\alpha{\v}T\onesccn.
\end{equation}
In this setting, $\gamma$ is a probability to stay in ESCC for one step if initial distribution is $\v$. For given $\v$, this number increases as $T$ becomes closer to stochastic matrix. Now, the problem of choosing $\gamma$ comes down to the problem of choosing $\v$. The advantage of this approach is twofold. First, we still have all the flexibility because, depending on $\v$, the value of $\gamma$ may be  literally anything between zero and one. Second, we can use a probabilistic interpretation of $\v$ to make a reasonable choice. In this paper we consider three appealing choices of $\v$: 
\begin{enumerate}
\item $\tilde{\pi}_T$, the quasi-stationary distribution  of $T$,
\item the uniform vector $\u$, and
\item the normalized PageRank vector $\pi_T(c)/||\pi_T(c)||_1$. Note that in this case both $\v$ and $\gamma$ depend on $c$.
\end{enumerate}

First, let us take $\v=\tilde{\pi}_T$, the quasi-stationary distribution of $T$. The motivation for taking $\v=\tilde{\pi}_T$ is that $\tilde{\pi}_T$ weights the states according to their quasi-stationary probabilities, which captures the structure of $T$. 
With $v=\tilde{\pi}_T$, equation (\ref{c*}) becomes
\[\pi_{T}(c^*)\onesccn= \alpha\tilde{\pi}_T T\onesccn=\alpha\lambda_1.
\]
In this case, $\gamma=\lambda_1$  is the probability that the random walk stays in ESCC given that it did not leave this block for infinitely long time. Hence, $\lambda_1$  is a natural measure of proximity of $T$ to stochastic matrix.

If conditions of Proposition~\ref{prop1} are satisfied, then (\ref{ub}) and (\ref{lb})  hold, and thus the value of $c^*$ satisfying (\ref{c*}) must be in the interval $(c_1,c_2)$, where 
\[(1-c_1)/(1-p_1 c_1)=\lambda_1,\quad (1-c_2)/(1-\lambda_1 c_2)=\lambda_1.\]
It is easy to check that $c_1=(1-\lambda_1)/(1-\lambda_1p_1)$ and $c_2=1/(\lambda_1+1)$. Since $\lambda_1$ is very close to 1, it follows that $c$ is bounded from above by the number $c^*\le c_2$, where $c_2$ is only slightly larger than $1/2$! Numerical results for our two datasets are presented in Table~\ref{tab:c}. 

\begin{table}[htb]
\centerline{
\begin{tabular}{|c|c|r|r|}
\hline
$v$&$c$&INRIA&FrMathInfo\\
\hline \hline
$\tilde{\pi}_T$&$c_1$&0.0184&0.1571\\
&$c_2$&0.5001&0.5002\\
&$c^*$&.02&.16\\
\hline
$\u$&$c_3$&0.5062&0.5009\\
&$c_4$&0.9820&0.8051\\
&$c^*$&.604&.535\\
\hline
$\pi_T/||\pi_T||_1$&$1/(1+\lambda_1)$&0.5001&0.5002\\
&$1/(1+p_1)$&0.5062&0.5009\\
\hline
\end{tabular}}
\caption{Values of $c^*$ with bounds.}
\label{tab:c}
\end{table}

From the numerical results we can see that in case when $v=\tilde{\pi}_T$, we obtain $c^*$ close to zero. This however leads to ranking that takes into account only local information about the Web Graph. Specifically, the number of incoming links will play a dominant role in the PageRank value (see e.g. \cite{Fortunato2}). Furthermore, the interpretation of $\v=\tilde{\pi}_T$ also suggests that it is not the best choice because the `easily bored surfer' random walk that is used in PageRank computations never follows a quasi-stationary distribution. Indeed, with probability $(1-c)$, this random walk restarts itself from the uniform probability vector. Clearly, the intervals between subsequent restarting points are too short to reach a quasi-stationary regime.

Our second choice is the uniform vector $\v=\u$. In this case, (\ref{c*}) becomes
\[||\pi_{T}(c^*)||_1= \alpha\u T\onesccn=\alpha p_1.
\]
If the conditions of Proposition~\ref{prop1} hold then we again can use (\ref{ub}) and (\ref{lb}) to establish that $c^*\in (c_3,c_4)$, where
\[(1-c_3)/(1-p_1c_3)=p_1, \quad (1-c_4)/(1-\lambda_1c_4)=p_1.\]
The values of $c_3=1/(1+p_1)$, $c_4=(1-p_1)/(1-\lambda_1p_1)$, and $c^*$ for our datasets are given in Table~1.
As we see, in this case, we have obtained a higher upper bound. However, the values of $c^*$ are still much smaller than $0.85$. 

\begin{figure}[hbt]
              \centering {\epsfxsize=2.8in \epsfbox{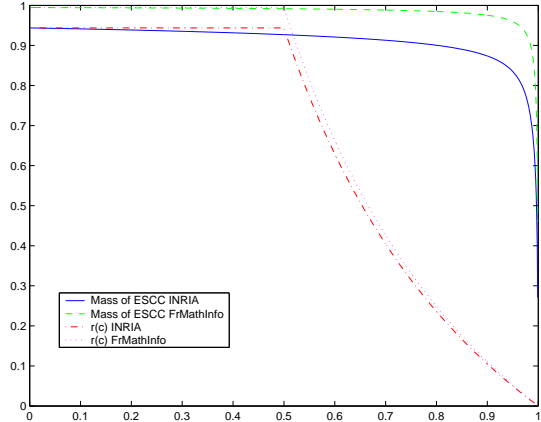}}
                \caption{\small The value $c^*$ with $v=\pi_{T}/||\pi_{T}||_1$ is the crossing point of $r(c)$ and $||\pi_{T}(c)||_1$.}
        \label{fig:plots_for_c}
\end{figure}

 Note that $\v=\tilde{\pi}_T$ implies $\gamma=\lambda_1$, which is a probability to stay in ESCC for one step after infinitely long time, and $\v=\u$ leads to $\gamma=p_1$, which is the probability to stay in ESCC for one step after starting afresh. Our third choice, the normalized PageRank vector $\v=\pi_T/||\pi_T||_1$,  is a symbiosis of the previous two cases. With this choice of $\v$, according to  (\ref{c*}), the value $c=c^*$ solves the equation
\begin{align*}
||\pi_{T}(c)||_1&= \frac{\alpha}{||\pi_{T}(c)||_1}\pi_{T}(c)T\onesccn\\
&= \frac{\alpha^2(1-c)}{||\pi_{T}(c)||_1}{\u}[I-cT]^{-1}T\onesccn,
\end{align*}
where the last equality follows from (\ref{PRESCC}).
Multiplying by $||\pi_{T}(c)||_1$,  we obtain
\begin{align*}
||\pi_{T}(c)||_1^2&=\alpha^2(1-c)\u\frac{1}{c}\,cT[I-cT]^{-1}\onesccn\\
&=\alpha^2(1-c)\u\frac{1}{c}\,cT[I-cT]^{-1}\onesccn\\
&=\alpha^2(1-c)\u\frac{1}{c}\,\left([I-cT]^{-1}-I\right)\onesccn\\
&=
\frac{\alpha}{c}\,||\pi_{T}(c)||_1-\frac{(1-c)\alpha^2}{c}.
\end{align*}
Solving the quadratic equation for $||\pi_{T}(c)||_1$, we get
\[||\pi_{T}(c)||_1=r(c)=\left\{\begin{array}{ll}\alpha&\mbox{if }c\le 1/2,\\
\frac{\alpha(1-c)}{c}&\mbox{if }c> 1/2.\end{array}\right.\]
Hence, the value $c^*$ solving (\ref{c*}) corresponds to the point where the graphs of $||\pi_{T}(c)||_1$ and $r(c)$ cross each other. First, note that there is only one such point on (0,1). Furthermore, since $||\pi_{T}(c)||_1$ decreases very slowly unless $c$ is close to one, and $r(c)$ starts decreasing relatively fast for $c>1/2$, one can expect that $c^*$ is only slightly larger than $1/2$. This is illustrated in Figure~\ref{fig:plots_for_c}, where we depict $||\pi_{T}(c)||_1$ and $r(c)$ for INRIA and FrMathInfo datasets. 

Under conditions of Proposition~\ref{prop1}, we may use (\ref{ub}) and (\ref{lb}) to deduce that $r(c)$ first crosses the line $\alpha(1-c)/(1-\lambda_1c)$, then $||\pi_{T}(c)||_1$, and then $\alpha(1-c)/(1-p_1c)$. Thus, we yield
\[\frac{1}{1+\lambda_1}<c^*< \frac{1}{1+p_1}.\]
Since both $\lambda_1$ and $p_1$ are close to 1, this clearly indicates that $c$ should be chosen close to 1/2. The values for lower and upper bounds for $c^*$ are given in Table~1. Since these bounds are tight we did not compute $c^*$ explicitly.

To summarize, our results indicate that with $c=0.85$, the ESCC component does not receive a fair share of the PageRank mass. Remarkably, in order to satisfy any of the three intuitive criteria of fairness presented above, the value of $c$ should be drastically reduced. In particular, the value $c=1/2$ looks like a well justified choice. 

In future, it would be interesting to design and analyze other criteria for choosing the most `fair' value of $c$. However, given the outcome of our studies, we foresee that any criterion based on the PageRank mass of ESCC will lead to similar results.

\section{Conclusions}
\label{sec:conclusions}

The choice of the PageRank damping factor is not evident. The old motivation for the value $c=0.85$ was a compromise between the true reflection of the Web structure and numerical efficiency. However, the Markov random walk on the Web Graph does not reflect the importance of the pages because it absorbs in dead ends. Thus, the damping factor is needed not only for speeding up the computations but also for establishing a fair ranking of pages. 

In this paper, we proposed new criteria for choosing the damping factor, based on the ergodic structure and probability flows. Our approach leads to the conclusion that the value $c=0.85$ is too high, and in fact the damping factor should be chosen close to 1/2. 

As we already mentioned before, the value $c=1/2$ was used in \cite{PRcitations} to find gems in scientific citations. This choice was justified intuitively by stating that researchers may check references in cited papers but on average they hardly go deeper than two levels, which results in probability 1/2 of `giving up'.  Nowadays, when search engines work really fast, this argument also applies to  Web search. Indeed, it is easier for the user to refine a query and receive a proper page in fraction of seconds than to look for this page by clicking on hyperlinks. Therefore, we may assume that a surfer searching for a page, on average, does not go deeper than two clicks.

Even if our statement that $c$ should be 1/2, might be received with a healthy skepticism, we hope to 
have convinced the reader that the study of ergodic structure of the Web helps in choosing the value 
of the damping factor, and in improving link-based ranking criteria in general. We believe that future research in this 
direction will yield new reasoning for a well grounded choice of the ranking criteria and help 
to discover new fascinating properties of the Web Graph.

\section*{Acknowledgments}

This work is supported by EGIDE ECO-NET grant no. 10191XC and by NWO Meervoud grant no.~632.002.401.
We also would like to thank Danil Nemirovsky for the collection of the Web Graph data.

\section*{Appendix: A Singular Perturbation Lemma}

\begin{lemma}
\label{lm:SPMC}
Let $G(\eps)=P+\eps C$ be a transition matrix of perturbed Markov chain.

The perturbed Markov chain is assumed to be ergodic for sufficiently small 
$\eps$ different from zero. And let the unperturbed Markov chain $(\eps =0)$
have $m$ ergodic classes. Namely, the transition matrix $P$ can be written
in the form

$$
P=\left[ \begin{array}{cccc} 
Q_1 &        & 0   & 0 \\
    & \ddots &     &   \\
0   &        & Q_m & 0 \\
R_1 & \cdots & R_m & T 
\end{array} \right] \in R^{n\times n}.
$$
Then, the stationary distribution of the perturbed Markov chain has a limit
$$
\lim_{\eps \to 0} \pi(\eps) =[\nu_1\mu_1 \ \cdots \ \nu_m\mu_m \ 0],
$$
where zeros correspond to the set of transient states in the unperturbed
Markov chain, $\mu_i$ is a stationary distribution of the unperturbed Markov
chain corresponding to the $i$-th ergodic set, and $\nu_i$ is the $i$-th
element of the aggregated stationary distribution vector that can be found
by solution
$$
\nu D = \nu, \quad \nu \one =1,
$$
where $D=MCQ$ is the generator of the aggregated Markov chain and
$$
M=\left[ \begin{array}{cccc} 
\mu_1 &        & 0   & 0 \\
    & \ddots &     &   \\
0   &        & \mu_m & 0  
\end{array} \right] \in R^{m\times n}.
$$
$$
Q=
\left[ \begin{array}{cccc} 
\one &        & 0    \\
    & \ddots &      \\
0   &        & \one  \\
\phi_1 & \cdots & \phi_m 
\end{array} \right] \in R^{n\times m}.
$$  
with $\phi_i=[I-T]^{-1}R_i\one$.
\end{lemma}
The proof of this lemma can be found in \cite{A99,KT93,YZ05}.

\end{document}